\documentclass[10pt,reqno]{amsart}

\usepackage{a4wide}
\usepackage{amssymb}
\usepackage{amsfonts}
\usepackage{amsmath}
\input xy
\xyoption{arrow} \xyoption{matrix}

\date{}

\newtheorem{proposition}{Proposition}[section]
\newtheorem{theorem}[proposition]{Theorem}
\newtheorem{lemma}[proposition]{Lemma}

\newtheorem{corollary}[proposition]{Corollary}

\def\Supp{\mathop\mathrm{Supp}}

\def\Hom{{\rm Hom}}
\def\der{\partial }

\def\nFM0{{\nu }_{F,M_0}}
\def\nFN0{{\nu }_{F,N_0}}
\def\nGN0{{\nu }_{G,N_0}}

\def\N0{ {\bf N}_0 }

\def\t{\otimes}
\def\g{\gamma}

\def\ra{\rightarrow}

\def\Xpm{X^{\pm }}

\def\Z{\mathbb{Z}}

\def\l1{{\lambda}_1}

\def\a{\alpha}
\def\a0{ {\alpha }_0}
\def\a1{ {\alpha }_1}

\def\l{\lambda}

\def\nFGM0{{\nu }_{F,G,M_0}}

\def\nFN0{{\nu}_{F,N_0}}

\def\sm{{\sigma}^m}

\def\sm1{{\sigma}^{-1}}

\def\smtp1{{\sigma}^{-t+1}}

\def\S1{S^{-1}}

\def\Xpm1{X^{\pm 1}_1}

\def\sPM1{{\sigma }^{\pm 1}}
\def\sMP1{{\sigma }^{\mp 1 }}

\def\d{\delta}

\def\di{{\rm d.ind}}

\def\L{\Lambda}

\def\Ytm1{Y^{t-1}}
\def\Yim1{Y^{i-1}}

\def\Aut{{\rm Aut}}

\def\dim{{\rm dim }}

\def\SL2Z{ {\rm SL}_2({\bf Z}) }

\def\Gp1{ G^{1 , 1 } }
\def\P11{ P^{-1 , 1 } }
\def\Pp1{ P^{1 , 1 } }

\def\Supp{{\rm Supp}}

\def\nCLsr{{}^\nu\kern-2pt {\cal L}^{\sigma , \rho  }}
\def\nP{{}^\nu \kern-2pt P}
\def\nL{{}^\nu\kern-2pt L}
\def\nLL{{}^\nu\kern-2pt \Lambda}
\def\nPsr{{}^\nu\kern-2pt P^{\sigma , \rho  }}
\def\nLsr{{}^\nu\kern-2pt L^{\sigma , \rho  }}
\def\nuCL{{}^\nu\kern-2pt  {\cal L}}
\def\nCLsr{{}^\nu\kern-2pt {\cal L}^{\sigma , \rho  }}
\def\nCL1m{{}^\nu\kern-2pt {\cal L}^{-1 , 1  }}
\def\x1nu{x^\frac{1}{\nu}}
\def\xm1nu{x^{-\frac{1}{\nu}}}

\def\ra{\rightarrow }

\def\CB{{\cal B}}

\def\nAM0{{\nu }_{{\cal A},M_0}}
\def\nAN0{{\nu }_{{\cal A},N_0}}

\def\End{ {\rm End }}

\def\SL{{\rm SL}}

\def\Ext{{\rm Ext}}
\def\Hom{{\rm Hom}}

\def\di!{\frac{\der^i}{i!}}
\def\dik!{\frac{\der^k_i}{k!}}

\def\N{\mathbb{N}}

\def\0{\overline{0}}
\def\1{\overline{1}}

\def\Ln1{\L_{n,\overline{1}}}

\def\a1{a_{\overline{1}}}

\def\S{\Sigma}

\def\vn1{\overrightarrow{n-1}}

\def\mJ{\mathbb{J}}
\def\mI{\mathbb{I}}

\def\bM{\overline{M}}

\def\K1{{\rm K}_1}

\def\hmI1{\widehat{\mI_1}}
\def\tmI1{\widetilde{\mI_1}}
\def\tmJ1{\widetilde{\mJ_1}}
\def\hB1{\widehat{B_1}}
\def\hCB1{\widehat{\CB_1}}

\begin{document}

\author{V. V.  \  Bavula, V.  Bekkert  \and V. Futorny}
\title[Indecomposable weight  modules]{Indecomposable generalized weight  modules over the algebra of polynomial
integro-differential operators}
\address{Department of Pure Mathematics, University of Sheffield, Hicks Building, Sheffield S3 7RH, UK}
\email{v.bavula@sheffield.ac.uk}
\address{Departamento de Matem\'atica, ICEx, Universidade Federal de Minas
Gerais, Av.  Ant\^onio Carlos, 6627, CP 702, CEP 30123-970, Belo
Horizonte-MG, Brasil} \email{bekkert@mat.ufmg.br}
\address{Instituto de Matem\'atica e Estat\'{\i}stica,
Universidade de S\~ao Paulo, Caixa Postal 66281, S\~ao Paulo, CEP
05315-970, Brasil} \email{futorny@ime.usp.br}



\begin{abstract}

For the algebra $\mI_1= K\langle x, \frac{d}{dx}, \int \rangle$ of
polynomial integro-differential operators over a field $K$ of
characteristic zero, a classification of  indecomposable, generalized weight $\mI_1$-modules  of finite length is
given. Each such module is an infinite dimensional uniserial module. Ext-groups are found between indecomposable generalized weight  modules, it is proven that they are finite dimensional vector spaces. 




\end{abstract}

\maketitle


\section{Introduction}

Throughout, ring means an associative ring with $1$; module means
a left module;
 $\N :=\{0, 1, \ldots \}$ is the set of natural numbers; $\N_+ :=\{ 1,2,  \ldots \}$ and $\Z_{\leq 0} :=-\N $;   $K$ is a
field of characteristic zero and  $K^*$ is its group of units;
$P_1:= K[x]$ is a polynomial algebra in one variable $x$ over $K$;
$\der:=\frac{d}{d x}$; $\End_K(P_1)$ is the algebra of all
$K$-linear maps from $P_1$ to $P_1$,  and $\Aut_K(P_1)$ is its
group of units (i.e. the
 group of all the invertible linear maps from $P_1$ to $P_1$); the
subalgebras  $A_1:= K \langle x , \der \rangle$ and
 $\mI_1:=K\langle x,\der ,  \int\rangle $ of $\End_K(P_1)$
  are called the (first) {\em Weyl algebra} and the {\em algebra  of polynomial
integro-differential operators} respectively  where $\int: P_1\ra
P_1$, $ p\mapsto \int p \, dx$, is the  {\em integration},  i.e.
$\int : x^n \mapsto \frac{x^{n+1}}{n+1}$ for all $n\in \N$.  The
algebra $\mI_1$ is neither left nor right Noetherian and not a
domain. Moreover, it contains infinite direct sums of nonzero left
and right ideals, \cite{algintdif}.

$\noindent $

In Section \ref{CCSI1MOD}, a classification of indecomposable, generalized weight
$\mI_1$-modules of finite length is given (Theorem \ref{6Apr14}). A similar classification is given in \cite{Bav-Bek}
for the generalized Weyl algebras where a completely different approach was taken. Properties of the algebras $\mI_n:=\mI_1^{\t n}$ of polynomial integro-differential operators
 in arbitrary many variables are studied in \cite{algintdif} and \cite{indtif-bimod}. The groups $\Aut_{K-{\rm alg}}(\mI_n)$ are found in \cite{intdifaut}. The simple $\mI_1$-modules are
 classified in \cite{algintdifline}.

\

\noindent{\bf Acknowledgment}
The first author is grateful to the  University of S\~ao Paulo for
 hospitality during his visit and to Fapesp for financial support (processo    2013/24392-5). The third first author is
supported in part by the  CNPq  (301320/2013-6), by the 
Fapesp  (2014/09310-5).


\section{Classification of indecomposable, generalized weight  $\mI_1$-modules of finite length}\label{CCSI1MOD}


In this section, a classification  of indecomposable, generalized weight  $\mI_1$-modules of finite length is
given (Theorem \ref{6Apr14}).

As an abstract algebra, the algebra $\mI_1$  is generated by the elements $\der $, $H:=
\der x$ and $\int$ (since $x=\int H$) that satisfy the defining
relations, {\cite[Proposition 2.2]{algintdif}} (where $[a,b]:=ab-ba$): $$\der \int = 1,
\;\; [H, \int ] = \int, \;\; [H, \der ] =-\der , \;\; H(1-\int\der
) =(1-\int\der ) H = 1-\int\der .$$ The elements of the algebra
$\mI_1$, 
\begin{equation}\label{eijdef}
e_{ij}:=\int^i\der^j-\int^{i+1}\der^{j+1}, \;\; i,j\in \N ,
\end{equation}
satisfy the relations $e_{ij}e_{kl}=\d_{jk}e_{il}$ where $\d_{jk}$
is the Kronecker delta function. Notice that
$e_{ij}=\int^ie_{00}\der^j$. The matrices of the linear maps
$e_{ij}\in \End_K(K[x])$ with respect to the basis $\{ x^{[s]}:=
\frac{x^s}{s!}\}_{s\in \N}$ of the polynomial algebra $K[x]$  are
the elementary matrices, i.e.
$$ e_{ij}*x^{[s]}=\begin{cases}
x^{[i]}& \text{if }j=s,\\
0& \text{if }j\neq s.\\
\end{cases}$$
Let $E_{ij}\in \End_K(K[x])$ be the usual matrix units, i.e.
$E_{ij}*x^s= \d_{js}x^i$ for all $i,j,s\in \N$. Then
\begin{equation}\label{eijEij}
e_{ij}=\frac{j!}{i!}E_{ij},
\end{equation}
 $Ke_{ij}=KE_{ij}$, and
$F:=\bigoplus_{i,j\geq 0}Ke_{ij}= \bigoplus_{i,j\geq
0}KE_{ij}\simeq M_\infty (K)$, the algebra (without 1) of infinite
dimensional matrices.

\

{\bf $\Z$-grading on the algebra $\mI_1$ and the canonical form of
an integro-differential operator, \cite{algintdif}}. The algebra
$\mI_1=\bigoplus_{i\in \Z} \mI_{1, i}$ is a $\Z$-graded algebra
($\mI_{1, i} \mI_{1, j}\subseteq \mI_{1, i+j}$ for all $i,j\in
\Z$) where  
\begin{equation}\label{I1iZ}
\mI_{1, i} =\begin{cases}
D_1\int^i=\int^iD_1& \text{if } i>0,\\
D_1& \text{if }i=0,\\
\der^{|i|}D_1=D_1\der^{|i|}& \text{if }i<0,\\
\end{cases}
\end{equation}
 the algebra $D_1:= K[H]\bigoplus \bigoplus_{i\in \N} Ke_{ii}$ is
a commutative non-Noetherian subalgebra of $\mI_1$, $ He_{ii} =
e_{ii}H= (i+1)e_{ii}$  for $i\in \N $ (notice that
$\bigoplus_{i\in \N} Ke_{ii}$ is the direct  sum of non-zero
ideals of $D_1$); $(\int^iD_1)_{D_1}\simeq D_1$, $\int^id\mapsto
d$; ${}_{D_1}(D_1\der^i) \simeq D_1$, $d\der^i\mapsto d$,   for
all $i\geq 0$ since $\der^i\int^i=1$.
 Notice that the maps $\cdot\int^i : D_1\ra D_1\int^i$, $d\mapsto
d\int^i$,  and $\der^i \cdot : D_1\ra \der^iD_1$, $d\mapsto
\der^id$, have the same kernel $\bigoplus_{j=0}^{i-1}Ke_{jj}$.

Each element $a$ of the algebra $\mI_1$ is the unique finite sum
\begin{equation}\label{acan}
a=\sum_{i>0} a_{-i}\der^i+a_0+\sum_{i>0}\int^ia_i +\sum_{i,j\in
\N} \l_{ij} e_{ij}
\end{equation}
where $a_k\in K[H]$ and $\l_{ij}\in K$. This is the {\em canonical
form} of the polynomial integro-differential operator
\cite{algintdif}.

$\noindent $

Let $v_i:=\begin{cases}
\int^i& \text{if }i>0,\\
1& \text{if }i=0,\\
\der^{|i|}& \text{if }i<0.\\
\end{cases}$

Then $\mI_{1,i}=D_1v_i= v_iD_1$ and an element $a\in \mI_1$ is the
unique  finite  sum 
\begin{equation}\label{acan1}
a=\sum_{i\in \Z} b_iv_i +\sum_{i,j\in \N} \l_{ij} e_{ij}
\end{equation}
where $b_i\in K[H]$ and $\l_{ij}\in K$. So, the set $\{ H^j\der^i,
H^j, \int^iH^j, e_{st}\, | \, i\geq 1; j,s,t\geq 0\}$ is a
$K$-basis for the algebra $\mI_1$. The multiplication in the
algebra $\mI_1$ is given by the rule:
$$ \int H = (H-1) \int , \;\; H\der = \der (H-1), \;\; \int e_{ij}
= e_{i+1, j},$$ $$ e_{ij}\int= e_{i,j-1}, \;\; \der e_{ij}=
e_{i-1, j},\;\; e_{ij} \der = \der e_{i, j+1},$$
$$ He_{ii} = e_{ii}H= (i+1)e_{ii}, \;\; i\in \N, $$
where $e_{-1,j}:=0$ and $e_{i,-1}:=0$.  $\noindent $

 The algebra
$\mI_1$ has the only proper ideal $F=\bigoplus_{i,j\in \N}Ke_{ij}
\simeq M_\infty (K)$ and
 $F^2= F$. The factor algebra $\mI_1/F$ is canonically isomorphic to
the skew Laurent polynomial algebra $B_1:= K[H][\der, \der^{-1} ;
\tau ]$, $\tau (H) = H+1$, via $\der \mapsto \der$, $ \int\mapsto
\der^{-1}$, $H\mapsto H$ (where $\der^{\pm 1}\alpha = \tau^{\pm
1}(\alpha ) \der^{\pm 1}$ for all elements $\alpha \in K[H]$). The
algebra $B_1$ is canonically isomorphic to the (left and right)
localization $A_{1, \der }$ of the Weyl algebra $A_1$ at the
powers of the element $\der$ (notice that $x= \der^{-1} H$).

An $\mI_1$-module $M$ is called a {\em weight} module if $M=\oplus_{\l \in K}M_\l$ where $M_\l := \{ m\in M\, | \, Hm=\l m\}$.
An $\mI_1$-module $M$ is called a {\em generalized weight} module if $M=\oplus_{\l \in K}M^\l$ where $M^\l := \{ m\in M\, | \, (H-\l )^nm=0$ for some $n=n(m)\}$.
The set $\Supp (M):=\{\l \in K\, | \, M^\l\neq 0\}$ is called the {\em support} of the generalized weight module $M$. For all $\l\in K$ and $n\geq 1$,
$$\der^n M^\l\subseteq M^{\l-n}\;\; {\rm and}\;\; \int^n M^\l\subseteq M^{\l+n}.$$
Let $0\ra N\ra M\ra L\ra 0$ be a short exact sequence of $\mI_1$-modules. Then $M$ is a generalized weight module iff so are the modules $N$ and $L$, and in this case
$$\Supp (M)= \Supp (N)\cup \Supp (L).$$

For each $\mI_1$-module $M$, there is a short exact sequence of $\mI_1$-modules
\begin{equation}
 0\ra FM\ra M\ra \bM:=M/FM\ra 0 \label{FMM}
\end{equation}

 where

 (i) $F\cdot FM=FM$, and

 (ii) $F\cdot \bM=0$,

 and the properties (i) and (ii) determine the short exact sequence
(\ref{FMM}) uniquely, i.e. if $0\ra M_1\ra M\ra M_2\ra 0$ is
a short exact sequence of $\mI_1$-modules such that $FM_1=M_1$ and $FM_2=0$ then $M_1\cong FM$
and $M_2\simeq \bM$.

Notice that
\begin{equation}
 FM\simeq K[x]^{I}, \label{FMM1}
\end{equation}

\noindent i.e. the $\mI_1$-module $FM$ is isomorphic to the direct sum of $I$ copies of the simple weight
$\mI_1$-module $K[x]$. Clearly, $\bM$ is a $B_1$-module.

{\bf The indecomposable $\mI_1$-modules $M(n,\lambda)$}.
For $\lambda\in K$ and  a natural number $n\geq 1$, consider the $B_1$-module
\begin{equation}
 M(n,\lambda):= B_1\t_{K[H]} K[H)/(H-\lambda)^{n}. \label{FMM2}
\end{equation}

\noindent Clearly,
\begin{equation}
 M(n,\lambda)\simeq B_1/B_1 (H-\lambda)^{n}\simeq\mI_1/(F+\mI_1 (H-\lambda)^{n}). \label{FMM3}
\end{equation}

\noindent The $\mI_1$-module/$B_1$-module $M(n,\lambda)$ is a generalized weight module with $\Supp M(n,\lambda)=\lambda+\Z$,
\begin{equation}
 M(n,\lambda)=\bigoplus_{i\in \Z}M(n,\lambda)^{\lambda +i} \,\,\,{\rm and}\,\, \dim \, M(n,\lambda)^{\lambda+i}=n\;\; {\rm for\; all}\;\; i\in \Z . \label{FMM4}
\end{equation}

For an algebra $A$, we denote by $A-{\rm Mod}$ its module category. The next proposition describes the set of indecomposable, generalized weight $\mI_1$-modules of finite length $M$ with $FM=0$.

\begin{proposition} \label{a7Apr14} \begin{enumerate}
\item  $M(n,\lambda)$ is an indecomposable, generalized weight $\mI_1$-module of finite length $n$.
\item $M(n,\lambda)\cong M(m,\mu)$ if and only if $n=m$ and $\lambda -\mu \in \Z$.
\item Let $M$ be a generalized weight $B_1$-module of length $n$ (i.e. let $M$ be a generalized weight $\mI_1$-module such that $FM=0$, by (\ref{FMM})).
 Then $M$ is indecomposable if and only if
                              $M\simeq M(n,\lambda)$ for some $\lambda\in K$.
\end{enumerate}
\end{proposition}

{\it Proof}. 1. Since $(B_1)_{K[H]}=\oplus_{i\in\Z}\der^{i}K[H]$ is a free right $K[H]$-module, the functor
$$B_1\t_{K[H]}-\;: K[H]-{\rm Mod}\ra B_1-{\rm Mod},\;\; N\mapsto B_1\t_{K[H]}N,$$
is an exact functor. The $K[H]$-module $K[H]/(H-\l )^n$ is an indecomposable, hence the $B_1$-module $M(n,\l )$ is indecomposable and  generalized weight  of length $n$.

2. $(\Rightarrow )$ Suppose that $\mI_1$-modules $M(n,\l )$ and $M(m,\mu )$ are isomorphic. Then $\Supp (M(n,\l))=\Supp (M(n,\l))$, i.e. $\l +\Z =\mu +\Z$, i.e. $\l -\mu \in \Z$. Then $n=m$, by (\ref{FMM4}).

$(\Leftarrow )$ Suppose that $k:=\l -\mu \in \Z$ and $n=m$. We may assume that $k\geq 1$.  Using the  equality $(H-\l )^n \der^k =\der^k (H-\l -k)^n=\der^k(H-\mu )^n$, we see that the $B_1$-homomorphism
$$M(n,\l )=B_1/B_1(H-\l )^n \ra M(n,\mu )=B_1/B_1(H-\mu )^n, \;\;1+B_1(H-\l )^n\mapsto \der^k+B_1(H-\mu )^n,$$ is an isomorphism with the inverse given by the rule $1+B_1(H-\mu )^n\mapsto \der^{-k}+B_1(H-\l )^n$.

3. $(\Leftarrow )$ This implication follows from statement 2.

$(\Rightarrow )$ Each indecomposable, generalized weight  $B_1$-module $M$ is of the type $B_1\t_{K[H]}N$ for an indecomposable $K[H]$-module $N$ of length $n$. Notice that  $N\simeq K[H]/(H-\l )^n$ for some $\l\in K$. Therefore, $M\simeq M(n, \l)$. $\Box$

\begin{lemma}\label{c7Apr14}
Let $M$ be an  indecomposable, generalized weight $\mI_1$-module. Then $\Supp (M)\subseteq \l +\Z$ for some $\l\in K$.
\end{lemma}

{\it Proof}. Let $M= \oplus_{\mu\in\Supp (M)}M^{\mu} $  be a generalized weight $\mI_1$-module. Then
$$M=\bigoplus_{\mu +\Z\in \Supp (M)/\Z} M_{\mu +\Z }$$ is a direct sum of
 $\mI_1$-submodules $M_{\mu +\Z }:=\oplus_{i\in\Z}M^{\mu +i}$ where $\Supp (M) /\Z$ is the image of the support $\Supp (M)$ under the abelian group epimorphism $ K\ra K/\Z$, $ \g\mapsto \g+\Z$.
 The $\mI_1$-module $M$ is indecomposable, hence $M=M_{\l+\Z}$ for some $\l\in K$, i.e. $\Supp (M)\subseteq \l +\Z$.
 $\Box$

The next lemma describes the set of indecomposable, generalized weight $\mI_1$-modules $M$ with $FM=M$.
\begin{lemma}\label{b7Apr14}
Let $M$ be an indecomposable, generalized weight $\mI_1$-modules $M$. Then the following statements are equivalent.
\begin{enumerate}
 \item $FM=M$.
 \item $M\simeq K[x]$.
 \item $\Supp (M) \subseteq \N$.
 \end{enumerate}
\end{lemma}

{\it Proof}. $(1)\Rightarrow (2):$ If $FM=M$ then $M\simeq K[x]^{(I)}$ for some set $I$ necessarily with $|I|=1$ since $M$ is indecomposable, i.e. $M\simeq K[x]$.

 $(2)\Rightarrow (3):$ $\Supp (K[x])=\{ 1,2,\ldots \}\subseteq \N$.

  $(3)\Rightarrow (1):$ Suppose that $\Supp (M)\subseteq \N$. Using the short exact sequence of $\mI_1$-modules $0\ra FM\ra M\ra \bM :=M/FM\ra 0$ we see that $\Supp (M) = \Supp (FM)\cup \Supp (\bM )$. Since $\Supp (FM)=\Supp (K[x]^{(I)})=\{ 1,2,\ldots \}$ and $\Supp (\bM )$ is an abelian group, we must have $\bM =0$ (since $\Supp (M) \subseteq \N$), i.e. $M=FM$.
$\Box$

The following result is a key step in obtaining a classification of indecomposable, generalized weight $\mI_1$-modules of finite length. 
\begin{theorem}\label{8Apr14}
Let $M$ be a generalized weight $\mI_1$-module of finite length. Then the short exact sequence (\ref{FMM}) splits.
\end{theorem}

{\it Proof}. We can assume that $FM\neq 0$ and $\bM \neq 0$. It is obvious that $FM\simeq K[x]^s$ for some $s\geq 1$ and the $B_1$-module $\bM \simeq \bigoplus_{i=1}^t M(n_i, \l_i )$ for some $n_i\geq 1$, $\l_i\in K$ and $t\geq 1$. It  suffices  to show that

\begin{equation}\label{ExtMK}
\Ext^1_{\mI_1}(M(n,\l ), K[x])=0
\end{equation}
for all $n\geq 1$ and $\l\in K$. If $\l\in \Z$ we can assume that $\l =0$, by Proposition \ref{a7Apr14}.(2).

(i) $F(H-\l )^n=F$: The equality follows from the equalities $e_{ij}(H-\l )^n=e_{ij}(j+1-\l )$ and the choice of $\l$.

(ii) $M(n,\l )=\mI_1/\mI_1(H-\l )^n$: By (i), $\mI_1(H-\l )^n \supseteq F(H-\l )^n = F$. Hence,
$$M(n,\l ) = \mI_1/(F+\mI_1(H-\l )^n)=\mI_1/\mI_1(H-\l )^n.$$

(iii) {\em The equality (\ref{ExtMK}) holds}: Let $M= M(n,\l )$. By (ii), the short exact sequence of $\mI_1$-modules
\begin{equation}\label{Mnlres}
0\ra \mI_1(H-\l )^n \ra \mI_1\ra M\ra 0
\end{equation}
is a projective resolution of the $\mI_1$-module $M$ since the map
$$\cdot (H-\l )^n : \mI_1\ra \mI_1(H-\l )^n,\;\; a\mapsto a(H-\l )^n,$$ is an isomorphism of $\mI_1$-modules, by the choice of  $\l$.  Then
$$\Ext^1_{\mI_1}(M, K[x])\simeq Z^1/B^1$$ where
$Z^1=\Hom_{\mI_1}(\mI_1(H-\l )^n, K[x])\simeq K[x]$ and $B^1\simeq (H-\l )^nK[x]=K[x]$, by the choice of $\l $.
 Hence,  the equality (\ref{ExtMK}) holds. The proof of the theorem is complete.
$\Box$

The next theorem is a classification of the set of indecomposable, generalized weight $\mI_1$-modules of finite length.

\begin{theorem}\label{6Apr14} 
Each indecomposable, generalized weight $\mI_1$-module of finite length is isomorphic to one of the modules below:
\begin{enumerate}
\item $K[x]$,
\item $M(n, \l )$ where $n\geq 1$ and $\l \in \L $ where $\L$ is any fixed subset of $K$ such that the map $\L\ra (K/\Z)$, $\l\mapsto \l +\Z$, is a bijection.
\end{enumerate}
The $\mI_1$-modules above are pairwise non-isomorphic, indecomposable, generalized weight and of finite length.
\end{theorem}

{\it Proof}. The theorem follows from Theorem \ref{8Apr14}, Proposition \ref{a7Apr14} and Lemma \ref{b7Apr14}.
$\Box$

\begin{corollary}\label{f8Apr14}
Every indecomposable, generalized weight $\mI_1$-module is an uniserial module. 
\end{corollary}

{\it Proof}. The statement follows from Theorem \ref{6Apr14}.$\Box $


\

{\bf Homomorphisms and Ext-groups between indecomposables.}

\begin{proposition}\label{b8Apr14}

\begin{enumerate}
\item Let $M$ and $N$ be generalized weight $\mI_1$-modules such that $\Supp (M)\cap \Supp (N)=\emptyset$. Then $\Hom_{\mI_1}(M,N)=0$.
\item  $\Hom_{\mI_1}(M (n,\l ) , K[x])=0$.
\item  $\Hom_{\mI_1}(K[x], M (n,\l ) )=0$.
\item $\Hom_{\mI_1}(M(n,\l ) , M(m, \l ))\simeq \Hom_{K[H]}(K[H]/((H-\l )^n), (K[H]/((H-\l )^m))\simeq$\\ $ K[H]/((H-\l )^{\min (n,m)})$.
\end{enumerate}
\end{proposition}

{\it Proof}. 1. Statement 1 is obvious.

2. Statement 2 follows from the fact that $FM (n,\l )=0$ and $Fp=K[x]$ for all nonzero elements $p\in K[x]$ (since $K[x]$ is a simple $\mI_1$-module, $F$ is an ideal of the algebra $\mI_1$ such that $FK[x]=K[x]$).

3. Statement 3 follows from the fact that  $FK[x]=K[x]$ and $FM(n,\l ) =0$: $f(K[x])=f(FK[x])=Ff(K[x])=0$
for any $f\in \Hom_{\mI_1}(K[x], M (n,\l ) )$.

4. The first isomorphism is obvious. Then the second isomorphism follows.  $\Box $

\begin{proposition}\label{e8Apr14}
\begin{enumerate}
\item  $\Ext_{\mI_1}^1(K[x],K[x])=0$.
\item  $\Ext_{\mI_1}^1(M(n,\l ),K[x])=0$.
\item  $\Ext_{\mI_1}^1(K[x],M(n,\l ))=0$.
\item $\Ext_{\mI_1}^1(M(n,\l ), M(m, \mu ))=\begin{cases}
K& \text{if } \l -\mu \in \Z ,\\
0& \text{if } \l -\mu \not\in \Z .\\
\end{cases}$
\end{enumerate}
\end{proposition}

{\it Proof}. 1.  Let $0\ra K[x]\ra N\ra K[x]\ra 0$ be a s.e.s. of $\mI_1$-modules. Then $FN=N$ (since $FK[x]=K[x]$), and so $N$ is an epimorphic image of the semisimple $\mI_1$-module $F\oplus F$. Hence, $N\simeq K[x]\oplus K[x]$ (since ${}_{\mI_1}F\simeq K[x]^{(\N )}$). 

2. See (\ref{ExtMK}). 

3. Let $0\ra M= M(n,\l ) \ra L\ra K[x]\ra 0$ be a s.e.s. of $\mI_1$-modules. Since $FM=0$, we have $FL=FK[x]\simeq K[x]$ is a submodule of $L$ such that $FL\cap M=0$ (since otherwise $FL\subseteq M$ by  simplicity of the $\mI_1$-module $FL\simeq K[x]$, and so $0\neq K[x]\simeq FL= F^2L\subseteq FM=0$, a contradiction). Then $FL\oplus M\subseteq L$. Furthermore, $FL\oplus M = L$ since $l_{\mI_1}(FL\oplus M)= l_{\mI_1}(L)$. This means that the s.e.s. splits. 

4. Let $0\ra M_1\ra  M\ra M_1\ra 0$ be a s.e.s. of generalized weight $\mI_1$-modules. If $\Supp (M_1)\cap \Supp (M_2) = \emptyset$, it splits. In particular, $\Ext^1_{\mI_1}(M(n,\l ), M(m, \mu ))=0$ if $\l -\mu \not\in \Z$. If $\l -\mu \in \Z$ we can assume that $\l =\mu$  (sine $M(m, \l) \simeq M(m,\mu)$). Using (\ref{Mnlres}), where we assume that $\l =0$ if $\l \in \Z$, we see that $\Ext^1_{\mI_1}(M(n,\l ), M(m, \l  ))\simeq M(m , \l)/(H-\l ) M(m,\l)\simeq K$.  $\Box$

Since the left global dimension of the algebra $\mI_1$ is 1, \cite{gldim-intdif}, 
Proposition \ref{b8Apr14} and Proposition \ref{e8Apr14} describe all the Ext-groups between indecomposable, generalized weight $\mI_1$-modules. This is also obvious from the proofs of the propositions.

\small{







}

\end{document}